\documentclass[12pt]{article}
\usepackage{amsmath,amssymb}
\newtheorem{theorem}{Theorem}[section]
\newtheorem{lemma}{Lemma}[section]

\newcommand{\mor}{morphism}
\newcommand{\x}{\times}
\renewcommand{\a}{\alpha}

\newcommand{\bs}{\bigskip}
\newcommand{\dt}{\cdot}
\newcommand{\C}{\mathbb C}

\newcommand{\e}{\varepsilon}

\renewcommand{\i}{\infty}
\renewcommand{\l}{\lambda}
\newcommand{\mb}{\mbox}
\newcommand{\ot}{\otimes}
\newcommand{\p}{\pi}
\newcommand{\sig}{\sigma}
\renewcommand{\t}{\tau}

\newcommand{{\z}}{\mathbb Z}
\newcommand{\R}{\mathbb R}
\newcommand{\frp}{free probability analogue }
\newcommand{\was}{Wasserstein }
\newcommand{\ib}{probability }
\newcommand{\res}{respectively }
\newcommand{\ts}{\tilde\Sigma}
\begin{document}
\setlength{\baselineskip}{18pt}

\begin{center}
{\large\bf A Free Probability Analogue of the}\\
{\large\bf Wasserstein Metric on the Trace-State Space}

\bs\bs {\sc Philippe Biane and Dan Voiculescu}

\bs{\it Incomplete preliminary version, June 1, 2000}
\end{center}

\bs\bs\bs\begin{quote}
{\bf Abstract.} We define a \frp of the  Wasserstein metric,
which extends the classical one. In dimension one, we prove that the square of
the Wasserstein distance to the semi-circle distribution is majorized by a
modified free entropy quantity.\end{quote}

\setcounter{section}{-1}
\section{Introduction}

The \was distance between two \ib distributions $\mu,\nu$ on ${\R}^n$
is given by
\[
W(\mu,\nu)=\inf_{\p\,\in\,\Pi(\mu,\nu)} (\int |x-y|^2\, d\p (x,y))^{\frac 12}
\]
where $\Pi(\mu,\nu)$ denotes the \ib measures on ${\R}^n\x{\R}^n$ with
marginals $\mu$ and $\nu$.

Following the usual free \ib recipe we shall replace the set of \ib measures
by the trace-state space of a $C^*$-algebra and take marginals with respect to
a free product. In this note we begin the study of the ensuing free \was
metric.

An inequality of M.~Talagrand ([7],[10]) relates the \was distance from a
Gaussian distribution and relative entropy. In the one-variable case we prove
a related free \ib result to this inequality, where the semicircle law
replaces the Gauss law and the logarithmic energy plays the role of entropy.
Note that in the case of $n$-tuples of commuting selfadjoint variables the
classical and the free \was distances are equal.

In the context of non-commutative geometry, there is a different
noncommutative extension, due to A.~Connes [5], of the related
Monge-Kantorowitz metric. The Monge-Kantorowitz metric is a $p\!=\!1$,
$p$-\was metric, but the definition which is extended is the dual definition
based on Lipschitz functions, and the extension involves Fredholm-modules or
derivations (recent work is surveyed in [9]).

\section{The free \was metric}

\subsection{The distance on $n$-tuples of variables}
We will work in the framework of tracial $C^*$-\ib spaces $(M,\tau)$, where
$M$ is a unital $C^*$-algebra and $\tau$ is a trace state. The simplest is to
define the metric at the level of noncommutative random variables.
If $(X_1,\dots ,X_n)$ and
$(Y_1,\dots ,Y_n)$ are two $n$-tuples of noncommutative random variables in
tracial
$C^*$-\ib spaces $(M_1,\tau_1)$ and $(M_2,\tau_2)$, we define
\[
W_p((X_1,\dots ,X_n),(Y_1,\dots ,Y_n))
\]
as the infimum of
\[
\|(|X'_j-Y'_j|_p)_{1\leq j\leq n}\|_p
\]
over $2n$-tuples $(X'_1,\dots ,X'_n,Y'_1,\dots ,Y'_n)$ of
  noncommutative random variables in some
  tracial $C^*$-\ib space $(M_3,\tau_3)$ such that the $n$-tuples
$(X'_1,\dots ,X'_n),(X_1,\dots ,X_n)$ and respectively
$(Y'_1,\dots ,Y'_n),(Y_1,\dots ,Y_n)$ have the same $*$-distributions.
Here $|\,\dt\,|_p$ is the $p$-norm in a  tracial $C^*$-\ib space,
while $\|\,\dt\,\|_p$ is the $p$-norm on ${\R}^n$. Like in the classical case,
if $p=2$ we call $W_p$ the free \was metric and we will also use the notation
$W$ for $W_2$. We shall refer to $W_p$ as the free $p$-\was metric. Note also
that if
\[
X_j = D_j+iE_j, \qquad Y_j=F_j+iG_j
\]
where $D_j,E_j,F_j,G_j$ are self-adjoint, then
\[
W((X_1,\dots ,X_n),(Y_1,\dots ,Y_n)) = W((D_1,\dots ,D_n,E_1,\dots ,E_n),
(F_1,\dots ,F_n,G_1,\dots ,G_n))
\]
Note also that $W_p((X_1,\dots ,X_n),(Y_1,\dots ,Y_n))$ depends only on the
$*$-distributions of $(X_1,\dots ,X_n)$ and $(Y_1,\dots ,Y_n)$. If we consider
$n$-tuples with the same $*$-distribution as equivalent; then $W_p$ will be a
distance between equivalence classes of
$n$-tuples.

\subsection{The distance on trace states}
We pass now to trace-state spaces $T\!S(A)$, where $A$ is a unital
$C^*$-algebra. We will assume $A$ is finitely generated  and we will assume
such a generator $(a_1,\dots ,a_n)$ has been specified. The $p$-\was metric
on $T\!S(A)$ is given by
\[
W_p(\tau',\tau'')=W_p((a'_1,\dots ,a'_n),(a''_1,\dots ,a''_n))
\]
where $\t',\t''\in TS(A)$ and $(a'_1,\dots ,a'_n)$ and $(a''_1,\dots ,a''_n)$
denote the variables defined by $(a_1,\dots ,a_n)$  in $(A,\t')$ and
respectively $(A,\t'')$.

This definition can be rephrased using free products. If $A_1,A_2$ are unital
$C^*$-algebras, we denote by $\sig_j:A_j\to A_1*A_2$ the canonical injection
of $A_j$ into the full free product $C^*$-algebra (this presumes amalgamation
over ${\C}1$). If $\t_j\in TS(A_j)$, $(1\leq j\leq 2)$ we define
\[
T\!S(A_1*A_2;\t_1,\t_2) =\{\t\in T\!S(A_1*A_2)\mid
\t\circ\sig_j=\t_j, \ j\!=\!1,2\} \ .
\]
Remark that $\t_1*\t_2\in T\!S(A_1*A_2;\t_1,\t_2)$.

It is easy to see that
\[
W_p(\t',\t'')=\inf\{\|(|\sig_1(a_j)-\sig_2(a_j)|_{p,\t})_{{}_{1\leq j\leq
n}}\|_{{}_p}
\mid \tau\in T\!S(A*A;\t',\t'')\}
\]
where $|\,\dt\,|_{p,\tau}$ denotes the $p$-norm in $L^p(A;\tau)$.

Remark also that the distance on $n$-tuples of variables can be obtained from
the definition for trace-states. Assume for simplicity $X_j=X^*_j$, $Y_j=Y^*_j$
and $R\geq\|X_j\|$, $R\geq\|Y_j\|$, $1\leq j\leq n$. Let then
$A=(C[-R,R])^{*n}$ (the free product of $n$ copies) and
$\sig_k(a)=~a_k$, where $a$ is the identical function in $C[-R,R]$. Let
$\rho_j:A\to M_j$, $j\!=\!1,2$ be the\linebreak $*$-homo\mor s such that
$\rho_1(a_k)=X_k$, $\rho_2(a_k)=Y_k$ where the $X_k$'s are in
$(M_1,\tau_1)$ and the $Y_k$'s in $(M_2,\t_2)$. Then
\[
W_p(\t',\t'')=W_p((X_1,\dots ,X_n),(Y_1,\dots ,Y_n))
\]
where $\t'=\t_1\circ\rho_1$, \ $\t''=\t_2\circ\rho_2$.

\bs\noindent{\bf 1.3 \ Theorem.} $W_p$ {\it is a metric.}

\bs {\bf Proof.} To check that $W_p$ is a metric on the set of equivalence
classes of $n$-tuples of variables or equivalently on a trace-state space
$ST(A)$ like in 1.2, the nontrivial assertion is the triangle inequality.
Indeed that $W_p((X_1,\dots ,X_n),(Y_1,\dots ,Y_n))\!=\!0\Leftrightarrow
(X_1,\dots ,X_n),(Y_1,\dots ,Y_n)$ have the same $*$-distribution or
\[
W_p(\t',\t'')=0\Leftrightarrow \t'=\t''
\]
are easy to see. For the  triangle inequality it will suffice to prove it in
the context of 1.1.

Let $(X'_1,\dots ,X'_n,Y'_1,\dots ,Y'_n)$ in $(M_{12},\t_{12})$ and
$(Y''_1,\dots ,Y''_n,Z''_1,\dots ,Z''_n)$ in $(M_{23},\t_{23})$ be $2n$-tuples
in tracial $W^*$-\ib spaces such that
$(X'_1,\dots ,X'_n)\sim (X_1,\dots ,X_n)$,\linebreak
$(Y'_1,\dots ,Y'_n)\sim (Y''_1,\dots ,Y''_n)\sim (Y_1,\dots ,Y_n)$,
$(Z''_1,\dots ,Z''_n)\sim (Z_1,\dots ,Z_n)$, where $\sim$ means the $n$-tuples
have equal $*$-distribution. There is a trace-preserving auto\mor \
of\linebreak
$W^*(Y'_1,\dots ,Y'_n)$ and $W^*(Y''_1,\dots ,Y''_n)$ which identifies $Y'_j$
and $Y''_j$. Abusing notations we shall denote  by $M_2$ the von Neumann
subalgebras of $M_{12}$ and $M_{23}$ generated by $(Y'_1,\dots ,Y'_n)$ and
respectively $(Y''_1,\dots ,Y''_n)$ identified as above. Let $E'$ and $E''$ be
the conditional expectations of $M_{12}$ and \res $M_{23}$ onto $M_2$.

Let $(M_{123},E)=(M_{12},E')*_{M_2} (M_{23},E'')$ and
$\t_{123}=\t_2\circ E$ where $\t_2=\t_{12}|M_2=\t_{23}|M_2$
(see 3.8 in [14]). Further, with $\rho_{12}:M_{12}\to M_{123}$,
$\rho_{23}:M_{23}\to M_{123}$ denoting the canonical embeddings, let
$X'''_j=\rho_{12}(X'_j)$, $Z'''_j=Z''_j$. Then
$\rho_{12}(Y'_j)=\rho_{23}(Y''_j)$ implies
\[
|X'''_j-Z'''_j|_{p,\t_{123}}\leq
|X'''_j-\rho_{12}(Y'_j)|_{p,\t_{123}} +
|\rho_{23}(Y''_j)-Z'''_j|_{p,\t_{123}}
= |X'_j-Y'_j|_{p,\t_{12}} + |Y''_j-Z''_j|_{p,\t_{23}}
\]
which is precisely what we need to establish the triangle inequality
\[
W_p((X_1,\dots ,X_n),(Y_1,\dots ,Y_n))+
W_p((Y_1,\dots ,Y_n),(Z_1,\dots ,Z_n))\geq
W_p((X_1,\dots ,X_n),(Z_1,\dots ,Z_n)) \ .
\]
${}$\hfill$\square$

\bs Let us also record as a proposition some easy consequences of the
compacity of the trace-state space. The proof is left to the reader.

\bs\noindent{\bf 1.4 \  Proposition.} (a) {\it The infimum in the definition of
$W_p$ is attained (both in the $1.1$ and $1.2$ contexts).}

(b) {\it Let \ $\t^{(k)}_1,\t_1,\t_2^{(k)},\t_2\in TS(A)$ \ and assume
$\t_j^{(k)}$ converges weakly to $\t_j$ as $k\to \i$ $(j=1,2)$. Then}
\[
\liminf_{k\to\i} W_p(\t_1^{(k)},\t_2^{(k)})\geq W_p(\t_1,\t_2) \ .
\]

(c) {\it Let} $(X_1^{(k)},\dots ,X_n^{(k)}),(X_1,\dots ,X_n)$,
$(Y_1^{(k)},\dots ,Y_n^{(k)}),(Y_1,\dots ,Y_n)$ {\it be $n$-tuples of\linebreak
variables in tracial $C^*$-\ib spaces and assume that}
$\|X_j^{(k)}\| \leq R$,
$\|X_j\| \leq R$,\linebreak  $\|Y_j^{(k)}\| \leq R$,
$\|Y_j\|\leq R$, {\it and that $(X_1^{(k)},\dots ,X_n^{(k)})$,
$(Y_1^{(k)},\dots ,Y_n^{(k)})$ converge in $*$-distribution to
$(X_1,\dots ,X_n)$ 	and \res $(Y_1,\dots ,Y_n)$. Then}
\[
\liminf_{k\to\i} W_p
((X_1^{(k)},\dots ,X_n^{(k)}),(Y_1^{(k)},\dots ,Y_n^{(k)}))\geq
W_p((X_1,\dots ,X_n),(Y_1,\dots ,Y_n)) \ .
\]

\bs If $(X_1,\dots ,X_n)$ are commuting self-adjoint variables in a tracial
$C^*$-\ib space, then their distribution $\mu_{X_1,\dots ,X_n}$ is a compactly
supported \ib measure on ${\R}^n$.

\bs\noindent{\bf 1.5 \ Theorem.} {\it Let} $(X_1,\dots ,X_n)$ and 
$(Y_1,\dots ,Y_n)$
{\it be two $n$-tuples of commuting self-adjoint variables in tracial $C^*$-\ib
spaces. Then the free and classical \was distances are equal:}
\[
W((X_1,\dots ,X_n),(Y_1,\dots ,Y_n)) =
W(\mu_{X_1,\dots ,X_n},\mu_{Y_1,\dots ,Y_n}) \ .
\]

\bs{\bf Proof.} The left-hand side is $\leq$ the right-hand side, since the
classical \was distance can be defined the same way as  the free one, with the
only difference that the $2n$-tuples $(X'_1,\dots ,X'_n,Y'_1,\dots
,Y'_n)$ in the infimum are required to live in commutative tracial $C^*$-\ib
spaces. We therefore only need to prove $\geq$\,.

Let $(X'_1,\dots ,X'_n,Y'_1,\dots ,Y'_n)$ be a $2n$-tuple in the infimum
defining the free distance. Passing to the von Neumann algebra completion, we
may assume $(M_3,\t_3)$, where $X'_j,Y'_j$ live, is a $W^*$-\ib space with a
normal faithful trace state. Let $A=W^*(X'_1,\dots ,X'_n)\subset M_3$,
$B=W^*(Y'_1,\dots,Y'_n)\subset M_3$ and let $E_A$ be the canonical conditional
expectation onto $A$. Then the unital trace-preserving completely positive map
$\varphi=E_A|B:B\to A$ gives  rise to a state $\nu:A\ot B\to{\C}$, on a
commutative algebra, defined by
\[
\nu(a\ot b)=\t_3(a\varphi(b)) \ .
\]
The positivity of $\nu$,
\[
\t_3(\sum_{i,j}a_ia^*_j\varphi(b_ib^*_j))\geq 0 \ ,
\]
is easily inferred from the positivity of the matrix
$(\varphi(b_ib^*_j))_{i,j}$. Alternatively, probabilistically, $\nu$ is the \ib
measure on ${\R}^{2n}$ obtained by integrating  w.r.t.
$\mu_{X_1,\dots ,X_n}$ the kernel of \ib measures describing
\[
\varphi: L^{\i}({\R}^n,\mu_{Y_1\dots Y_n})\to
  L^{\i}({\R}^n,\mu_{X_1\dots X_n}) \ .
\]
Then
\begin{eqnarray*}
\sum_{1\leq j\leq n}\nu ((X'_j-Y'_j)^2) &=&
\sum_{1\leq j\leq n}\t_3 (X_j^{'2}+\varphi(Y_j^{'2})-X'_j\dt
\varphi(Y'_j)-Y'_j\varphi(X'_j))  \\  &=&
\sum_{1\leq j\leq n}\t_3 (X'_j+Y_j^{'2}-2E_A(X'_jY'_j))\\  &=&
\sum_{1\leq j\leq n}\t_3 ((X'_j-Y'_j)^2) \ .
\end{eqnarray*}
Since $A\ot B$ is commutative this proves the
theorem.\hfill$\square$

\bs
\section{Cost of transportation to the semicircle distribution}

\subsection{The complex quasilinear differential equation}
Let $X,S$ in $(M,\t)$ be self-adjoint and freely independent and assume $S$ is
(0,1) semicircular. The purpose of section 2 is to estimate $W(X,S)$.
We begin by  studying variables $X(t)=e^{-t/2}X+(1-e^{-t})^{\frac 12}S$
which have the same distribution as the variables in the free
Ornstein-Uhlenbeck process. For technical reasons, and without extra work, the
complex PDE will be derived under the more general assumption that $X$ is
unbounded self-adjoint affiliated with $M$ (see [1]).

If $Y$ is self-adjoint affiliated with $M$, we denote by $\mu_X$ its
distribution and by $G_{\mu_{{}_Y}}(z)$ or $G_Y(z)$ the Cauchy transform of
$\mu_Y$, which equals $\tau((zI-X)^{-1})$.

If $Y(r)=X+r^{\frac 12}S$, let \
${\tilde G}(r,z)=G_{Y(r)}(z)$ and $G(t,z)=G_{X(t)}(z)$,
Im $z>0$, $r\geq 0$, $t\geq 0$. Then ${\tilde G}$ satisfies the complex
Burgers equation (see [3],[12])
\[
\frac{\partial{\tilde G}}{\partial r} +
{\tilde G} \ \frac{\partial{\tilde G}}{\partial z} =0 \ .
\]
Like ${\tilde G}(t,z)$ also $G(t,z)$ is $C^1$ on $[0,\i)\x\{z\in{\C}\mid
{\mb{Im}} z>0\}$ and holomorphic in $z$ for fixed $t$.
Note that $X(t)=e^{-t/2}Y(e^t)$ and that $G_{\a Y}(z)=\a^{-1}G(\a^{-1}z)$.
It follows that $G(t,z)=e^{t/2}{\tilde G}(e^t,e^{t/2}z)$. The complex Burgers
equation then gives
\begin{equation}
\frac{\partial G}{\partial t} +(G-\frac z2)
\frac{\partial G}{\partial z} -\textstyle{\frac 12} G=0
\end{equation}
with initial data $G(0,z)=G_X(z)$.

\subsection{The transport equation} Here we shall assume that the distribution
of $X$ is of the form $P_{\l}*\mu$ where $P_{\l}$ is the Cauchy distribution
with density $\p^{-1}\l(\l^2+x^2)^{-1}$ $(\l>0)$ and $\mu$ has compact
support. Since  $P_{\l}*\mu=P_{\l}\boxplus\mu$ ([1]) this is equivalent to
replacing $X$ with $X+\l C$ where $X$ is bounded, $X$ and $C$ are free and $C$
has a  Cauchy distribution $P_1$. Note that
$\mu_{X+\l C+r^{\frac 12}S} =\mu_{X+r^{\frac 12}S}*P_{\l}$, \
$G_{X+\l C+r^{\frac 12}S}(z)=G_{X+r^{\frac 12}S}(z+i\l)$, etc.
Thus, if the distribution of $X$ is of the form $P_{\l}*\mu$ then the equation
(1) is satisfied on an extended domain
\[
\{(t,z)\in [0,\i)\x{\C}\mid {\mb{Im}} \ z> -e^{t/2}\l\} \ .
\]

Let $-\p^{-1}G(x,t)=q(x,t)+ip(x,t)$ where $x\in{\R}$. Then $p(\dt \, ,t)$ is
the density of $\mu_{X(t)}$ and is analytic. For fixed $t$ and $k\geq 0$ we
have
\[
\left| \frac{\partial^k}{\partial x^k} \ p(x,t)\right| =
O((1+|x|)^{-2-k}) \quad {\mb{and}}\quad
\left| \frac{\partial^k}{\partial x^k} \ q(x,t)\right| =
O((1+|x|)^{-1-k}) \ .
\]
Moreover these bounds are uniform for $t$ in a compact set.

Equation (1) gives
\begin{eqnarray}
q_t &=& \p(qq_x-pp_x)+ 2^{-1}(xq_x+q) \nonumber \\
p_t &=& \p(pq_x +qp_x)+2^{-1}(xp_x+p) \\
q &=& -Hp \nonumber
\end{eqnarray}
where $H$ denotes the Hilbert transform.

Since $p(x,t)>0$ we infer that
\[
f(a,t)=\int^a_{-\i} p(x,t)dx
\]
is a $C^{\i}$-diffeo\mor s $f(\dt \,t):{\R}\to (0,1)$ which transports
$\mu_{X(t)}$ to Lebesgue measure. Hence
$\varphi_{s,t}(\dt)=f^{-1}(f(\dt \,s),t)$ \ $(0< s <t)$ will be a
  $C^{\i}$-diffeo\mor \ ${\R}\to{\R}$, which transports  $\mu_{X(s)}$ to
$\mu_{X(t)}$. This is the  same as saying that $X(t)$ and
  $\varphi_{s,t}(X(s))$ have the same distribution.

It is easily seen that
\[
\frac{\partial}{\partial t} f^{-1}(y,t) \ = \
\frac{-(\frac{\partial}{\partial t} f)(f^{-1}(y,t),t)}{p(f^{-1}(y,t),t)} \ .
\]
Using (2) to compute $\frac{\partial}{\partial t} f$ we find
\[
\frac{\partial}{\partial t} f(a,t) =\int^a_{-\i}
(\p(pq)_x +2^{-1}(xp)_x)dx \ = \ \p(pq)(a,t)+2^{-1}ap(a,t) \ .
\]
Hence
\[
\frac{\partial}{\partial t} f^{-1}(y,t) \ = \
-\p q(f^{-1}(y,t),t)-2^{-1}f^{-1}(y,t) \ .
\]
For \ $y=f(x,s)$ \ we get the transport equation
\begin{equation}
\frac{\partial}{\partial t} \varphi_{s,t}(x) \ = \
\p\big(Hp(\dt \, , t))(\varphi_{s,t}(x))-2^{-1}\varphi_{s,t}(x)\big)
\end{equation}
with initial condition $\varphi_{s,s}(x)=x$.

By the $L^m$-continuity $(1<m<\i)$ results for the density
(see Corollary 2 in [\quad ]) applied to $\mu\boxplus\mu_{r^{\frac 12}S}$ as a
function of $r$, we infer after convolutions with Cauchy distributions the
continuity of
\[
(0,\i)\ni t\longrightarrow Hp(\dt \, ,t)\in L^m({\R})
\]
(the $L^m$-space w.r.t.~Lebesgue measure). The reader should keep these facts
in mind in computations where we shall use (3).

\setcounter{lemma}{2}
\begin{lemma} Assume $X$ has distribution $\mu *P_{\l}$, where $\mu$ has
compact support and let\linebreak $X(t)=e^{-t/2}X+(1+e^{-t})^{\frac 12}S$ with
$S$ \ $(0,1)$-semicircular and free from $X$. Let $g\in C^{\i}({\R})$ be such
that $\|g\|_{\i}<\i$, \ $\|g'\|_{\i}\leq 1$ and assume $g'$ has compact
support. Then
\[
(t-s)^2W(g(X(s)),g(X(t)))^2\leq \sup_{s\leq h\leq t} \
\int_{\mb{\rm{\footnotesize supp}} \ g'}
(\p Hp(\dt \, ,h)(x)-2^{-1}x)^2 p(x,h)dx \ .
\]
\end{lemma}

{\bf Proof.} We have
\begin{eqnarray*}
&& W(g(X(s)),g(X(t)))^2 \\
&&\qquad\leq \int_{\R} |g(x)-g(\varphi_{s,t}(x))|^2 p(x,s)dx \\
&&\qquad \leq \int_{\R} \left( \int^t_s g'(\varphi_{s,h}(x))
(\p Hp(\dt \, ,h)(\varphi_{sh}(x))-2^{-1}\varphi_{s,h}(x))dh\right)^2p(x,s)dx\\
&&\qquad\leq  (t\!-\!s)\int_{\R}\int^t_s
(g'(\varphi_{s,h}(x)))^2
(\p Hp(\dt \, ,h)(\varphi_{sh}(x))-2^{-1}\varphi_{s,h}(x))^2 dh \ p(x,s)dx\\
&&\qquad= (t\!-\!s)\int^t_s\left(\int_{\R}
(g'(\varphi_{s,h}(x)))^2
(\p Hp(\dt \, ,h)(\varphi_{s,h}(x))-2^{-1}\varphi_{s,h}(x))^2 dh \
p(x,s)dx\right)dh\\ &&\qquad =
  (t\!-\!s)\int^t_s \int_{\R} (g'(x))^2
(\p Hp(\dt \, ,h)(x)-2^{-1}x)^2 \ p(x,h)dxdh\\
&&\qquad\leq (t-s)^2\sup_{s\leq h\leq t} \
\int_{\mb{\rm{\footnotesize supp}} \ g'}
(\p Hp(\dt \, , h))(x)-2^{-1}x)^2 p(x,h)dx \ .
\end{eqnarray*}
${}$\hfill$\square$

\bs
{\large\bf 2.4.} Assume $X$ is bounded and the semicircular variable $S$ is
free w.r.t.~$X$. Then the distribution $\mu_{X(t)}$ of
$X(t)=e^{-t/2}X+(1-e^t)^{\frac 12}S$ has $L^{\i}$-density $p(\dt \, ,t)$
w.r.t.~Lebesgue measure (see any of the papers [1],[2],[3],[11],[12]).

\bs{\bf Lemma.} {\it Assume $X$ is bounded, $S$ is $(0,1)$ semicircular, $X$
and
$S$ are free and let $p(\dt \, ,t)$ be the density of $\mu_{X(t)}$, where
$X(t)=e^{-t/2}X+(1-e^t)^{\frac 12}S$. Then}
\[
(t-s)^{-2}W(X(s),X(t))^2\leq \sup_{s\leq h\leq t} \
\int(\p Hp(\dt \, ,h)(x)-2^{-1}x)^2 p(x,h)dx \ .
\]

\bs

{\bf Proof.} Let $C$ be a variable with Cauchy distribution and free
w.r.t~$\{X,S\}$. Let $g\in C^{\i}({\R})$ be such that $\|g'\|_{\i}\leq 1$,
$g(x)=x$ \ if \ $|x|\leq\|X\|+1$ and $g'(x)=0$ \ if \ $|x|\geq\|X\|+2$.
We shall apply Lemma 2.3 to $X+\l C$ in place of $X$. Let
\[
Z(t,\l)=e^{-t/2}(X+\l C)+(1-e^{-t})^{\frac 12}S=X(t)+e^{-t/2}\l C \ .
\]
Then $g(Z(t,\l))$ is an operator of norm $\leq\|X\|+2$ and converges in
distribution to $X(t)$. Moreover the distribution of $Z(t,\l)$ is given by
the density $P_{e^{-t/2}\l}* p(\dt \, ,t)$ and will be denoted by
$p(\dt \, ,t,\l)$. In view of the $L^m$-continuity of $p(\dt \, ,t)$
$(1<m<\i)$ ([12]) it is easy to see that
\begin{eqnarray*}
&& \limsup_{\l\downarrow 0} (\sup_{s\leq h\leq t} \
\int_{\mb{\rm{\footnotesize supp}} \ g'}
(\p Hp(\dt \, ,h,\l)(x)-2^{-1}x)^2 p(x,h,\l)dx  \\
&&\qquad \leq \sup_{s\leq h\leq t} \
\int (\p Hp(\dt \, ,h)(x)-2^{-1}x)^2 p(x,h)dx \ .
\end{eqnarray*}
${}$\hfill$\square$

\bs
{\large\bf 2.5.} From now on we return to the context of bounded variables
$X$. If the distribution of $X$ is Lebesgue absolutely continuous and has
density $p$ which is $L^3$, then $\frac 12 {\mathcal J}(X)=\p Hp(X)$ where
${\mathcal J}(X)$ is the conjugate variable (a.k.a.~free Brownian gradient,
a.k.a.~noncommutative Hilbert transform) (see [13]) and
\[
\Phi(X)=\t({\mathcal J}(X)^2) =
4\p^2\int (Hp(x))^2p(x)dx=
\textstyle{\frac 43} \p^2\!\!\int p^3(x)dx
\]
is the free Fisher information (see [11],[13] up to  different
normalizations). The quantity occurring in Lemma 2.4,
\[
I(X)=4\!\int (\p Hp(x)-2^{-1}x)^2 \ p(x)dx =\t(({\mathcal J}(X)-X)^2) =
\Phi(X)-2+\t(X^2) \ ,
\]
is a generalization of the free Fisher information for Ornstein-Uhlenbeck
processes (see [4]). The inequality in Lemma 2.4 can also be written
\begin{equation}
4(t-s)^{-2} W(X(s),X(t))^2\leq \sup_{s\leq h\leq t}I(X(h)) \ .
\end{equation}

\setcounter{subsection}{5}
\subsection{The free entropy} The free entropy of $X$ with distribution
$\mu=\mu_X$ is
\[
\chi(X) = \iint\log |s-t|d\mu(s)d\mu(t) +
\textstyle{\frac 34}+\textstyle{\frac 12}\log(2\p)
\]
(see [11],[13] up to  different constants) and we have
\[
\chi(\a X)=\chi(X)+\log|\a| \quad {\mb{and}}\quad
\lim_{\e\downarrow 0} \e^{-1}(\chi (X+\e^{\frac 12}S)-\chi(X))=
2^{-1}\Phi (X) \ .
\]

The quantity we shall use in estimating the distance to the semicircle
distribution is a modified free entropy adapted to the free
Ornstein-Uhlenbeck process ([4]):
\begin{eqnarray*}
{\tilde\Sigma}(X) & = &-\chi(X) + \chi(S) +\textstyle{\frac 12}
\t(X^2)-\textstyle{\frac 12} \\
&=& \textstyle{\frac 12}
\t(X^2)-\displaystyle{\iint} d\mu(s)d\mu(t)\log |
s-t|-\textstyle{\frac 34} \ .
\end{eqnarray*}
We have \ $\displaystyle{\lim_{t\to\i}}{\tilde\Sigma}(X(t))=0$ \ and
\begin{eqnarray*}
\frac{d}{dt} {\tilde\Sigma}(X(t)) & = &
\frac{d}{dt}\left( \frac t2-
\chi(X+(e^t-1)^{\frac 12} S) + \textstyle{\frac 12} e^{-t}\t(X^2)+
\textstyle{\frac 12}(1-e^{-t})\right) \\
&=& 2^{-1} \left( 1-e^t\Phi(X+(e^t-1)^{\frac 12}S)
-e^{-t}\t(X^2)+e^{-t}\right) \\
&=& 2^{-1}(1-\Phi(X(t))-\t(X(t)^2)+1) = -2^{-1} I(X) \- .
\end{eqnarray*}

Note also that in [4] using the logarithmic Sobolev inequality for $\chi$
(Prop.~7.9 in [13]), it is shown that
\begin{equation}
{\tilde\Sigma}(X(t))\leq 2^{-1} I(X(t))
\end{equation}
which is a logarithmic Sobolev inequality for the Ornstein-Uhlenbeck process.

\setcounter{lemma}{6}
\begin{lemma} Assume $X,Y$ are bounded and self-adjoint, then if $t>0$ we have
\[
\limsup_{\e\to 0} |\e|^{-1}
|W(Y,X(t+\e))-W(Y,X(t))| \leq 2^{-1}(I(X(t)))^{\frac 12}
\]
\end{lemma}

{\bf Proof.} By the triangle inequality for $W$, we have
\[
|W(Y,X(t+\e))-W(Y,X(t))| \leq W(X(t),X(t+\e)) \ .
\]
The lemma then follows from (4) and the continuity of $I(X(h))$ $(h>0)$,
which is a consequence of the continuity of $\Phi(X(h))$
(Corollary 2 in [12]).\hfill$\square$

\bs We now have all ingredients to get an estimate for $W(X,S)$ which is
similar in the free context to an inequality of Talagrand in the classical
setting ([7],[10]).

\setcounter{theorem}{7}
\begin{theorem} $W(X,S)^2\leq{\tilde\Sigma}(X)$.
\end{theorem}

{\bf Proof.} Because of the semicircular maximum for $\chi$ we have
$\chi(X)\leq\chi(S) +2^{-1}\log (\t(X^2))$ so that
${\tilde\Sigma}(X)\geq 2^{-1}(\t(X^2)-(1+\log\t(X^2)))\geq 0$. Thus it will
suffice to prove that $W(X,S)-({\tilde\Sigma}(X))^{\frac 12}\leq 0$.

By Lemma 2.7, the inequality (5) and the formula for the derivative of
$\ts(X(t))$, we have for $t>0$,
\begin{eqnarray*}
&& \liminf_{\e\to 0} \e^{-1} (W(X(t+\e),S)-
(\ts (X(t+\e)))^{\frac 12}-W(X(t),S)+(\ts(X(t)))^{\frac 12}) \\
&&\qquad \geq -2^{-1}(I(X(t)))^{\frac 12}+2^{-2}I(X(t))
(\ts(X(t)))^{-\frac 12} \\
&&\qquad \geq 2^{-1}(I(X(t)))^{\frac 12}+2^{-2+1}
I(X(t)) (I(X(t)))^{-\frac 12} \ = \ 0 \ .
\end{eqnarray*}
Hence $W(X(t),S)-(\ts(X(t)))^{\frac 12}$ is an increasing function
and we have
\[
\lim_{t\to\i} (W(X(t),S)-(\ts(X(t)))^{\frac 12}) \ = \ 0
\]
because of the semicircular maximum and lower semicontinuity of $\chi$.
It follows that
\[
W(X(t),S)-(\ts(X(t)))^{\frac 12}\ \leq \  0
\]
if $t>0$.  To get the inequality for $t=0$, remark that $X(t)$ is
norm-continuous so that $W(X(t),S)$ tends to $W(X,S)$ as $t\to 0$.
On the other hand, by lower semicontinuity of $\chi$,
\[
\liminf_{t\downarrow 0} (-(\ts(X(t)))^{\frac 12})\ \geq \
-(\ts(X))^{\frac 12} \ .
\]
${}$\hfill$\square$

\setcounter{subsection}{8}
\bs\subsection{Remark.} Because of the coincidence of the free and classical
\was distance for single self-adjoint variables, the preceding theorem
can also be written in terms of \ib measures for the classical distance.
Let $\mu$ be a compactly supported \ib measure on ${\R}$ and
$\sig$ a (0,1)-semicircle distribution. Then we have
\[
(W(\mu,\sig))^2\leq\textstyle{\frac 12}\int x^2d\mu(x)-
{\displaystyle{\iint}} d\mu(s)d\mu(t)\log |s-t|-\textstyle{\frac 34} \ .
\]

\bs\bs\begin{quote}{\bf Acknowledgment.} This research was conducted by Dan
Voiculescu for the Clay Mathematics Institute. He was also supported in part by
National Science Foundation grant DMS95--00308.\end{quote}

\bs\bs\begin{center}{\bf References}\end{center}
\begin{description}
\item{[1]} Bercovici, H., Voiculescu, D.  Free convolution of
measures with unbounded support. {\it Indiana Univ. Math. J.} {\bf 42} (1993),
no.~3, 733--773.

\item{[2]} Biane, P. Processes with  free increments.
{\it Math. Z.} {\bf 227} (1998), 143--174.

\item{[3]}  Biane, P. On the free convolution with a semicircular distribution.
{\it Indiana Univ. Math. J.} {\bf 46} (1997), 705--717.

\item{[4]} Biane, P., Speicher R. Free diffusions, free entropy and free
Fisher information. Preprint DMA, \ ENS 99--33 (1999).

\item{[5]} Connes, A. Compact metric spaces, Fredholm modules and
hyperfiniteness. {\it Ergodic Theory and Dynamical Systems} {\bf 9} (1989),
207--220.

\item{[6]} Evans, L.C. Partial differential equations and
Monge-Kantorovich mass transfer. Preprint.

\item{[7]} Otto, F., Villani, C. Generalization of an inequality by
Talagrand, and links with the logarithmic Sobolev inequality.
  Preprint DMA, \ ENS 99--23 (1999).

\item{[8]} Rachev, S.T. {\it Probability Metrics and the Stability of
Stochastic Models}, Wiley series in \ib and mathematical statistics (1991).

\item{[9]} Rieffel, M.A. Metrics on state space. Preprint (1999).

\item{[10]} Talagrand, M. Transportation cost for Gaussian and other product
measures. {\it Geom. Funct. Anal.} {\bf 6} (1996), 587--600.

\item{[11]} Voiculescu, D. The analogues of entropy and of Fisher's
information measure in free \ib theory, I. {\it Commun. Math. Phys.} {\bf 155}
(1993), 71--92.

\item{[12]} Voiculescu, D. The derivative of order 1/2 of a free convolution
by a semicircle distribution. {\it Indiana Univ. Math. J.} {\bf 46} (1997),
697--703.

\item{[13]} Voiculescu, D. The analogues of entropy and of Fisher's
information measure in free \ib theory, V: noncommutative Hilbert transforms.
{\it Invent. Math.} {\bf 132} (1998), 189--227.

\item{[14]} Voiculescu, D., Dykema, K., Nica, A. {\it Free Random Variables.}
CRM Monograph Series no.~1, American Mathematical Society, Providence, RI
(1992).
\end{description}

\bs\bs
\mbox{}\hfill
\begin{tabular}{r}
Philippe Biane\\
{\sc cnrs, dma,} \'Ecole Normnale Sup\'erieure\\
45 Rue d'Ulm\\
75005 Paris, France\\
{\sc biane@dmi.ens.fr}\bigskip\\
Dan Voiculescu\\
Department of Mathematics\\
University of California\\
Berkeley, California 94720--3840 USA\\
{\tt dvv@math.berkeley.edu}
\end{tabular}

\end{document}